\documentclass{gtart}
\usepackage{amsmath,amscd,amsfonts}
\input gtmonout
\volumenumber{1}
\volumeyear{1998}
\volumename{The Epstein birthday schrift}
\pagenumbers{317}{334}
\papernumber{15}
\received{15 November 1997}
\published{26 October 1998}

\newtheorem{theorem}{Theorem}
\newtheorem{lemma}{Lemma}[section]
\newtheorem{corollary}[lemma]{Corollary}
\def\func#1{\text{\rm #1}\,}
\def\limfunc#1{\text{\rm #1}\,}

\theoremstyle{definition}

\begin{document}

\title{On the continuity of bending}

\author{Christos Kourouniotis}

\address{Department of Mathematics \\
University of Crete \\
Iraklio, Crete, Greece}
\email{chrisk@math.uch.gr}

\begin{abstract}
We examine the dependence of the deformation obtained by bending
quasi-Fuchsian structures on the bending lamination. We show that when we
consider bending quasi-Fuchsian structures on a closed surface, the
conditions obtained by Epstein and Marden to relate weak convergence of
arbitrary laminations to the convergence of bending cocycles are not
necessary. Bending may not be continuous on the set of all measured
laminations. However we show that if we restrict our attention to
laminations with non negative real and imaginary parts then the deformation
depends continuously on the lamination.
\end{abstract}

\primaryclass{30F40}\secondaryclass{32G15}

\keywords{Kleinian groups, quasi-Fuchsian groups, geodesic laminations}

\maketitle

The deformation of hyperbolic structures by bending along totally
geodesic submanifolds of codimension one was introduced by Thurston in his
lectures on \emph{The Geometry and Topology of 3--manifolds}. The geometric
and algebraic properties of the deformation were studied in \cite{K1} and 
\cite{JM}. Epstein and Marden \cite{EM} introduced the notion of a bending
cocycle and used it to describe bending a hyperbolic surface along a
measured geodesic lamination. The same notion was used in \cite{K2} to
extend bending to a holomorphic family of local biholomorphic homeomorphisms
of quasi-Fuchsian space $Q(S)$.

Epstein and Marden \cite{EM} give a careful analysis of the dependence of
the bending cocycle on the measured lamination. They consider the set of
measured laminations on $\mathcal{H}^2$ consisting of geodesics that
intersect a compact subset $K\subset \mathcal{H}^2$. This is a subset of the
space of measures on the space $G(K)$ of geodesics in $\mathcal{H}^2$
intersecting $K$, with the topology of weak convergence of measures. In this
topology, the bending cocycle does not depend continuously on the
lamination. One reason for this is the behaviour of the laminations near the
endpoints of the segment over which we evaluate the cocycle. For example,
consider the geodesic segment $[e^{i\theta },i]$ in $\mathcal{H}^2$, for
suitable $\theta $ in $[0,\pi /2]$, and the measured laminations $\mu _n$,
with weight $1$ on the geodesic $(1/n,n)$ and weight $-1$ on the geodesic $%
(-1/n,-n)$. Then $\{\mu _n\}$ converges weakly to the zero lamination, but
the cocycle of $\mu _n$ relative to $[e^{i\theta },i]$ is approximately a
hyperbolic isometry of translation length $1$. Epstein and Marden find
conditions under which a sequence of measured laminations gives a convergent
sequence of cocycles relative to a given pair of points.

In this article we show that when the lamination is invariant by a discrete
group and we only consider cocycles relative to points in the orbit of a
suitable point $x\in \mathcal{H}^2$, any sequence of measured laminations $%
\{\mu _n\}$ which converges weakly gives rise to cocycles which converge up
to conjugation. We show further that the same conjugating elements can be
used for the cocycles for $\mu _n$ corresponding to the different generators
of the group. Hence the laminations $\mu _n$ determine bending homomorphisms
which, after conjugation by suitable isometries, converge to the bending
homomorphism determined by $\mu _0$. This implies that the deformations
converge in $Q(S)$.

\begin{theorem}
\label{thm1}Let $S$ be a closed hyperbolic surface and $Q(S)$ its space of
quasi-Fuchsian structures. Let $\{\mu _n\}$ be a sequence of complex
measured geodesic laminations, converging weakly to a lamination $\mu _0$.
Then the bending deformations 
\[
B_{\mu _n}\co \mathcal{D}_{\mu _n}\rightarrow Q(S)
\]
converge to the deformation $B_{\mu _0}$, uniformly on compact subsets of $%
\mathcal{D=D}_{\mu _0}\cap \left( \bigcup_{m=1}^\infty \bigcap_{n=m}^\infty 
\mathcal{D}_{\mu _n}\right) $.
\end{theorem}

We also state an infinitesimal version of the Theorem.

\begin{theorem}
\label{thm2}Let $S$ be a closed hyperbolic surface and $Q(S)$ its space of
quasi-Fuchsian structures. Let $\{\mu _n\}$ be a sequence of complex
measured geodesic laminations, converging weakly to a lamination $\mu _0$.
Then the holomorphic bending vector fields $T_{\mu _n}$ on $Q(S)$ converge
to $T_{\mu _0}$, uniformly on compact subsets of $Q(S)$.
\end{theorem}

These results do not necessarily imply the continuous dependence of the
deformation on the bending lamination, because the space of measured
laminations is not first countable. If however we restrict our attention to
the subset of measured laminations with non negative real and imaginary
parts, then we can apply results in \cite{PH} to obtain the following
Theorem.

\begin{theorem}
\label{thm3}The mapping $\mathcal{ML}^{++}(S)\times Q(S)\rightarrow
T(Q(S))\co (\mu ,[\rho ])\mapsto T_{\mu} ([\rho ])$ is continuous, and
holomorphic in $[\rho ]$.
\end{theorem}

The proof of Theorem \ref{thm1} is based on the observation that, when the
lamination is invariant by a discrete group and we are considering cocycles
with respect to points $x$ and $g(x)$, for some $g$ in the group, the effect
of a lamination near the endpoints of the segment $[x,g(x)]$ is controlled
by its effect near $x$, provided that the lamination does not contain
geodesics very close to the geodesic carrying $[x,g(x)]$. This last
condition can be achieved by choosing $x$ to be a point not on the axis of a
conjugate of $g$ (see Corollary \ref{cor10}).

In Section \ref{sec1} we describe the space of measured laminations and we
recall the definition of bending. In the beginning of Section \ref{sec2} we
recall or modify certain results from \cite{EM} and \cite{K2} which provide
bounds for the effect of bending along nearby geodesics. Lemma \ref{lem9}
and the results following it examine the consequences of the above condition
on the choice of $x$.

The proof of Theorems \ref{thm1}, \ref{thm2} and \ref{thm3} is given in
Section \ref{sec3}. The laminations $\mu _n$ are replaced by finite
approximations. The main result is Lemma \ref{lemmain}, which gives the
basic estimate for the difference between the bending homomorphism of $\mu _0
$ and a conjugate of the bending homomorphism of $\mu _n$. Then a diagonal
argument is used to obtain the convergence of bending.

\section{The setting\label{sec1}}

We consider a closed surface $S$ of genus greater than 1. We fix a
hyperbolic structure on $S$, and let $\rho _0\co \pi _1(S)\rightarrow PSL(2,%
\mathbb{R)}$ be an injective homomorphism with discrete image $\Gamma _0=\rho
_0(\pi _1(S)$, such that $S$ is isometric to $\mathcal{H}^2/\Gamma _0$.

We consider the space $R$ of injective homomorphisms $\rho \co \Gamma
_0\rightarrow PSL(2,\mathbb{C)}$ obtained by conjugation with a quasiconformal
homeomorphism $\phi $ of $\widehat{\mathbb{C}}$: if $g\in \Gamma _0$, acting
on $\widehat{\mathbb{C}}$ as M\"{o}bius transformations, then $\rho (g)=\phi
\circ g\circ \phi ^{-1}$.

$PSL(2,\mathbb{C)}$ acts on the left on $R$ by inner automorphisms. The
quotient of $R$ by this action is the \emph{space} $Q(S)$ \emph{of
quasi-Fuchsian structures} on $S$, or \emph{quasi-Fuchsian space} of $S$. We
denote the equivalence class in $Q(S)$ of a homomorphism $\rho \in R$ by $%
[\rho ]$. Then $[\rho ]$ is a \emph{Fuchsian point} if there is a circle in $%
\widehat{\mathbb{C}}$ left invariant by $\rho (\Gamma _0)$, so that $\rho
(\Gamma _0)$ is conjugate to a Fuchsian group of the first kind. The subset
of Fuchsian points in $Q(S)$ is the \emph{Teichm\"{u}ller} space of $S$, $%
T(S)$.

We fix a point $[\rho ]\in Q(S)$, represented by the homomorphism $\rho
\co \Gamma _0\rightarrow$\break$PSL(2,\mathbb{C)}$ obtained by conjugation with the
quasiconformal homeomorphism $\phi \co \widehat{\mathbb{C}}\rightarrow \widehat{%
\mathbb{C}}$. We denote the image of $\rho $ by $\Gamma $. The limit set of $%
\Gamma _0$ is $\widehat{\mathbb{R}}$. Then $\phi (\widehat{\mathbb{R}})$ is the
limit set of $\Gamma $. If $\gamma $ is a geodesic in $\mathcal{H}^2$ with
endpoints $u,v\in \widehat{\mathbb{R}}$, we denote by $\phi _{*}(\gamma )$ the
geodesic in $\mathcal{H}^3$ with endpoints $\phi (u),\phi (v)$ in $\phi (%
\widehat{\mathbb{R}})$. In this way, geodesics on the surface $S\cong \mathcal{H%
}^2/\Gamma _0$ are associated to geodesics in the hyperbolic 3--manifold $%
\mathcal{H}^3/\Gamma $.

We want to study the deformation of quasi-Fuchsian structures by \emph{%
bending}, \cite{K1}, \cite{EM}, \cite{K2}. Bending is determined by a geodesic
lamination on $S$ with a complex valued transverse measure.

A measured geodesic lamination on $S$ lifts to a measured geodesic
lamination on $\mathcal{H}^2$. The space $G(\mathcal{H}^2)$ of unoriented
geodesics in $\mathcal{H}^2$ is homeomorphic to a M\"{o}bius strip without
boundary. Let $K$ be a compact subset of $\mathcal{H}^2$, projecting onto $%
\mathcal{H}^2/\Gamma _0$. The set $G(K)$ of geodesics in $\mathcal{H}^2$
intersecting $K$ is a compact metrizable space.

A measured geodesic lamination\emph{\ }on $\mathcal{H}^2$ determines a
complex valued Borel measure $\mu $ on $G(K)$, with the property that if $%
\gamma _1$ and $\gamma _2$ are distinct geodesics in the support of $\mu $,
then they are disjoint. The set of measured geodesic laminations on $S$ can
be considered as a subset of $\mathcal{M}(G(K))$, the set of complex valued
Borel measures on $G(K)$. The set $\mathcal{M}(G(K))$ has a norm, defined by 
\[
\left\| \mu \right\| =\sup \left\{ \left| \int f\mu \right| \text{, }f\text{
continuous complex valued function on }G(K)\text{,}\left| f\right| \leq
1\right\} 
\]
We shall use the weak* topology on $\mathcal{M}(G(K))$, with basis the sets
of the form 
\[
U(\mu ,\varepsilon ,f_1,\ldots ,f_m)=\left\{ \nu \in \mathcal{M}%
(G(K)):\left| \int f_i\mu -\int f_i\nu \right| <\varepsilon ,i=1,\ldots
,m\right\} 
\]
where $\mu \in \mathcal{M}(G(K))$, $f_i$, $i=1,\ldots ,m$ are continuous
functions on $G(K)$, and $\varepsilon $ is a positive number.

A measured geodesic lamination $\mu $ on $S$ is called \emph{finite} if it
is supported on a finite set of simple closed geodesics in $S$. Then, for
any compact subset $K$ of $\mathcal{H}^2$, the measure on $G(K)$ determined
by the lift of $\mu $ to $\mathcal{H}^2$ has finite support.

Given a finite measured geodesic lamination $\mu $ on $S$, we define bending
the quasi-Fuchsian structure $[\rho ]$ on $S$ as follows.

Let $g_1,\ldots ,g_k$ be a set of generators of $\Gamma _0$. Choose a point $%
x$ on $\mathcal{H}^2$ and, for each $g_j$, consider the geodesic segment $%
[x,g_j(x)]$. Let $\gamma _1,\ldots ,\gamma _m$ be the geodesics in the
support of $\mu $ intersecting $[x,g_j(x)]$, and let $z_1,\ldots ,z_m$ be
the corresponding measures. If $\gamma _1$ (or $\gamma _m$) go through $x$
(or $g_j(x)$ respectively), we replace $z_1$ (or $z_m$) by $\frac 12z_1$ (or 
$\frac 12z_m$).

If $\gamma $ is an oriented geodesic in $\mathcal{H}^3$ and $z\in \mathbb{C}$,
we denote by $A(\gamma ,z)$ the element of $PSL(2,\mathbb{C)}$ with axis $%
\gamma $ and complex displacement $z$. We will use the same notation for one
of the matrices in $SL(2,\mathbb{C)}$ corresponding to $A(\gamma ,z)$. In such
cases either the choice of the lift will not matter, or there will be an
obvious choice.

We orient the geodesics $\gamma _1,\ldots ,\gamma _m$ so that they cross the
segment $[x,g_j(x)]$ from right to left, and define the isometry 
\[
C_{t\mu }(x,g_j(x))=A(\phi _{*}(\gamma _1),tz_1)\cdots A(\phi _{*}(\gamma
_m),tz_m). 
\]
For each generator $g_j,$ $j=1,\ldots ,k$, define 
\[
\rho _{t\mu }(g_j)=C_{t\mu }(x,g_j(x))\text{ }\rho (g_j). 
\]
For $t$ in an open neighbourhood of $0$ in $\mathbb{C}$, the representation $%
[\rho _{t\mu }]$ is quasi-Fuchsian, \cite{K1}.

Any measured geodesic lamination $\mu $ on $S$ can be approximated by finite
laminations so that the corresponding bending deformations converge, \cite
{EM}, \cite{K2}. In this way, we obtain for any measured geodesic lamination
on $S$ a deformation $B_{\mu} $ defined on an open set $\mathcal{D}_{\mu}
\subset Q(S)\times \mathbb{C}$, 
\[
B_{\mu} \co \mathcal{D}_{\mu} \rightarrow Q(S)\co ([\rho ],t)\mapsto [\rho _{t\mu }]%
. 
\]
$B_{\mu} $ is a holomorphic mapping.

\section{The lemmata\label{sec2}}

In the vector space $\mathbb{C}^2$ we introduce the norm 
\[
\left\| (z_1,z_2)\right\| =\max \{|z_1|,|z_2|\}.
\]
A complex matrix $A=\left( 
\begin{array}{ll}
a & b \\ 
c & d
\end{array}
\right) $ acts on $\mathbb{C}^2$ and has norm 
\[
\left\| A\right\| =\max \left\{ |a|+|b|,|c|+|d|\right\}. 
\]
We will use this norm on $SL(2,\mathbb{C)}$.

\begin{lemma}[\cite{EM}, 3.3.1]
\label{lem1} Let $X$ be a set of matrices
in $SL(2,\mathbb{C)}$ and $c=(0,0,1)\in \mathcal{H}^3$. Then the following are
equivalent.
\begin{itemize}
\item[\rm i)] The closure of $X$ is compact.
\item[\rm ii)] There is a positive number $M$ such that if $A\in X$ then $||A||\leq M$.
\item[\rm iii)] There is a positive number $M$ such that if $A\in X$ then $||A||\leq M$
and $||A^{-1}||\leq M$.
\item[\rm iv)] There is a positive number $R$ such that if $A\in X$ then $d(c,A(c))\leq
R$.\endproof
\end{itemize}
\end{lemma}

Let $\Lambda $ be a maximal geodesic lamination on $S$, and $\psi
\co S\rightarrow \mathcal{H}^3/\Gamma $ the pleated surface representing the
lamination $\Lambda $ \cite{CEG}. Let $\tilde{\psi}\co \mathcal{H}%
^2\rightarrow \mathcal{H}^3$ be the lift of $\psi $.

\begin{lemma}[\cite{K2}, 2.5]
\label{lem1.1} Let $K$ be a compact disc of
radius $R$ about $c=(0,0,1)\in \mathcal{H}^3$, and $M$ a positive number.
There is a positive number $N$ with the following property. If $[x,y]$ is a
geodesic segment in $\mathcal{H}^2$ such that $\tilde{\psi}([x,y])\subset K$
and $\{\gamma _i,z_i\}$, $i=1,\ldots ,m$ is a finite measured lamination
with support contained in $\Lambda $, whose leaves all intersect $[x,y]$ and
are numbered in order from $x$ to $y$, and such that $\sum_{i=1}^m|\func{Re}%
z_i|<M$, then 
\[
\left\| A(\gamma _1,z_1)\cdots A(\gamma _m,z_m)\right\| \leq N.
\]
\vglue -.8cm
\endproof
\end{lemma}

\begin{lemma}[\cite{EM}, 3.4.1, \cite{K2}, 2.4]
\label{lem2} Let K be a
compact subset of $SL(2,\mathbb{C)}$, $M$ a positive number, and let $\gamma $
be the geodesic $(0,\infty )$. Then there is a positive number $N$ with the
following property. For any $B,C\in K$, and $z\in \mathbb{C}$ with $|z|\leq M$,
we have 
\[
\left\| BA(\gamma ,z)B^{-1}-CA(\gamma ,z)C^{-1}\right\| \leq N\left\|
B-C\right\| \left| z\right|.
\]
\vglue -.8cm 
\endproof
\end{lemma}

In order to examine the effect of bending along nearby geodesics, in Lemma 
\ref{lem4} and \ref{lem4.1}, we shall use the notion of a solid cylinder in
hyperbolic space. A \emph{solid cylinder }$C$ \emph{over a disk }$D$ in $%
\mathcal{H}^n$ is the union of all geodesics orthogonal to a $(n-1)$%
--dimensional hyperbolic disc $D$ in $\mathcal{H}^n$. The \emph{radius} of
the cylinder is the hyperbolic radius of the disc $D$. If $x$ is the centre
of $D$, we say that $C$ is a solid cylinder \emph{based} at $x$. The
boundary of $C$ at infinity consists of two discs $D_1$ and $D_2$ in $%
\partial \mathcal{H}^n$. We say that the solid cylinder $C$ is \emph{%
supported }by $D_1$ and $D_2$. The geodesic orthogonal to $D$ through its
centre is the \emph{core} of the solid cylinder $C$. We shall denote the
cylinder with core $\gamma $, basepoint $x\in \gamma $ and radius $r$ by $%
C(\gamma ,x,r)$.

\begin{lemma}[\cite{K2}, 2.6]
\label{lem3} Let $L$ be a compact set in $%
\mathcal{H}^3$. Then there exists a positive number $M$ with the following
property. If $D$ is a disc of radius $r$, contained in $L$, and $\alpha
,\beta $ are two geodesics contained in the solid cylinder over $D$, then
there is an element $A\in SL(2,\mathbb{C)}$ such that $A(\alpha )=\beta $ and $%
||A-I||\leq Mr$.\hfill 
\endproof%
\end{lemma}

If $C$ is a solid cylinder supported on the discs $D_1$ and $D_2$, with $%
D_1\cap D_2=\emptyset $, and $\gamma _1,\gamma _2$ are two geodesics, each
having one end point in $D_1$ and one in $D_2$, we say that $\gamma _1$ and $%
\gamma _2$ are \emph{concurrently oriented }in $C$ if their origins lie in
the same component of $D_1\cup D_2$.

\begin{lemma}
\label{lem4}Let $m$ be a positive number and $L$ a compact subset of $%
\mathcal{H}^3$. Then there are positive numbers $M_1$ and $M_2$ with the
following property. If $\gamma _1,\gamma _2$ are concurrently oriented
geodesics contained in a cylinder of radius $r$, based at a point in $L$,
and $z_1,z_2$ are complex numbers such that $|z_i|\leq m$, then there are
lifts of $A(\gamma _i,z_i)$ to $SL(2,\mathbb{C)}$ such that 
\[
\left\| A(\gamma _1,z_1)-A(\gamma _2,z_2)\right\| \leq M_1r\min
\{|z_1|,|z_2|\}+M_2|z_1-z_2|.
\]
\end{lemma}

\proof%
We assume that $|z_1|\leq |z_2|$. We have 
\[
\left\| A(\gamma _1,z_1)-A(\gamma _2,z_2)\right\| \leq \left\| A(\gamma
_1,z_1)-A(\gamma _2,z_1)\right\| +\left\| A(\gamma _2,z_1)-A(\gamma
_2,z_2)\right\|. 
\]
Let $B\in SL(2,\mathbb{C})$ be an element mapping the geodesic $(0,\infty )$ to 
$\gamma _2$, and mapping the point $c=(0,0,1)$ to a point in $L$. Then, by
Lemma \ref{lem1}, there is a constant $K_1$ depending only on $L$, such that 
$||B||\leq K_1$. By Lemma \ref{lem3} there is an element $C\in SL(2,\mathbb{C)}$
such that $C(\gamma _2)=\gamma _1$, and $||C-I||\leq K_2r$ for some constant 
$K_2$ depending only on $L$.

By Lemma \ref{lem2} there is a constant $K_3$ such that 
\[
\left\| A(\gamma _1,z_1)-A(\gamma _2,z_1)\right\| \leq K_3\left\|
CB-B\right\| |z_1|\leq K_1K_2K_3r|z_1|. 
\]
On the other hand, 
\[
\left\| A(\gamma _2,z_1)-A(\gamma _2,z_2)\right\| \leq \left\| B\right\|
\left\| A((0,\infty ),z_1-z_2)-I\right\| \left\| B^{-1}\right\| \left\|
A((0,\infty ),z_2)\right\| . 
\]
By Lemma \ref{lem1} and the fact that the entries of $A((0,\infty ),z_1-z_2)$
depend analytically on $z_1-z_2$, there is a constant $K_4$, depending on $L$
and $m$ such that 
\[
\left\| A(\gamma _2,z_1)-A(\gamma _{2,}z_2)\right\| \leq K_4|z_1-z_2|. 
\]
\vglue -.85cm
\endproof

\begin{lemma}[\cite{K2}, 2.7]
\label{lem4.1} Let $m$ be a positive number
and $L$ a compact subset of $\mathcal{H}^3$. Then there is a positive number 
$M$ with the following property. Let $C$ be a solid cylinder of radius $r$
based at a point in $L$. Let $\gamma _1,\ldots ,\gamma _k$ be geodesics in $%
C\,$and $z_1,\ldots ,z_k$ complex numbers with $\sum_{i=1}^k|\func{Re}%
(z_i)|\leq m$. Then 
\[
\left\| A(\gamma _1,z_1)\cdots A(\gamma _k,z_k)-A\left( \gamma
_1,\sum_{i=1}^kz_i\right) \right\| \leq Mr\sum_{i=1}^k|z_i|.
\]
\vglue -.85cm
\endproof
\end{lemma}

We want to show that if two geodesics on $S$ are sufficiently close, then
the corresponding geodesics in $\mathcal{H}^3/\Gamma $ will also be close,
(Lemma \ref{lem8}).

\begin{lemma}
\label{lem5}Let $K$ be a compact subset of $\mathcal{H}^2$, and $\phi
\co \partial \mathcal{H}^2\rightarrow \partial \mathcal{H}^3$ a homeomorphism
onto its image. Then there is a compact subset $L$ of $\mathcal{H}^3$ such
that if $\gamma $ is a geodesic of $\mathcal{H}^2$ intersecting $K$, then $%
\phi _{*}(\gamma )$ intersects $L$, i.e. $\phi _{*}(G(K))\subset G(L)$.
\end{lemma}

\proof%
We consider the Poincar\'{e} disk model of hyperbolic space. There, it is
clear that if $K$ is a compact subset of $B^2$, then there is a positive
number $m$ such that if $\gamma $ is a geodesic in $G(K)$ with end-points $%
u,v$, then $|u-v|\geq m$. Since $\phi ^{-1}$ is uniformly continuous, there
is a positive number $M$ such that $|\phi (u)-\phi (v)|\geq M$, and hence
there is a compact subset of $B^3$ intersecting $\phi _{*}(\gamma )$.\hfill 
\endproof%

\begin{lemma}[\cite{K2}, 2.2]
\label{lem6} Let $\varepsilon $ and $\eta $
be two positive numbers. Then there is a positive number $\delta $ with the
following property. If $D_1$ and $D_2$ are discs in $S^2$, with spherical
radius $\leq \delta $, and the spherical distance between $D_1$ and $D_2$ is 
$\geq \eta $, then the solid cylinder supported by $D_1$ and $D_2$ has
hyperbolic radius $r\leq \varepsilon $.\hfill 
\endproof%
\end{lemma}

\begin{lemma}
\label{lem7}Let $K$ be a compact subset of $B^n$, and $d$ a positive number.
Then there is a positive number $\delta $ with the following property. If $C$
is a solid cylinder in $B^n$, over a disc with radius $r\leq \delta $ and
centre at a point in $K$, then the spherical radius of each of the discs
supporting $C$ is $\leq d$.
\end{lemma}

\proof%
The radii of the supporting discs are given by continuous functions of the
core geodesic, the base point and the radius of the cylinder. For a fixed
base point, they tend to zero with the radius of the cylinder. The result
follows by compactness.\hfill%
\endproof%

\begin{lemma}
\label{lem8}Let $[\rho ]$ be a quasi-Fuchsian structure on $S$, $K$ a
compact subset of $\mathcal{H}^2$, and $L$ a compact subset of $\mathcal{H}^3
$ such that $\phi _{*}(G(K))\subset G(L)$. Let $r$ be a positive number.
Then there is a positive number $\delta $ with the following property. If $%
\gamma \in G(K)$, $x\in \gamma \cap K$ and $0\leq r_1\leq \delta $, then
there is some point $x^{\prime }\in L$ such that for any geodesic $\alpha $
contained in the solid cylinder $C(\gamma ,x,r_1)$, the geodesic $\phi
_{*}(\alpha )$ is contained in the solid cylinder $C(\phi _{*}(\gamma
),x^{\prime },r)\subset \mathcal{H}^3$.
\end{lemma}

\proof%
We work in the Poincar\'{e} disc model of the hyperbolic plane and space, $%
B^2$ and $B^3$. Since $L$ is a compact subset of $B^3$, there is a number $%
\eta _2>0$ such that if $u$ and $v$ are the endpoints of any geodesic in $%
B^3 $ intersecting $L$, then the spherical distance between $u$ and $v$ is $%
\geq \eta _2$. Then, by Lemma \ref{lem6}, there is a positive number $\delta
_2$, such that any solid cylinder with core a geodesic $\gamma \in G(L)$ and
supported on discs of spherical radius $\leq \delta _2$, has hyperbolic
radius $\leq r$.

Since $\phi \co S^1\rightarrow S^2$ is uniformly continuous, there is a
positive number $\delta _1$, such that any arc in $S^1$ of length $\leq
\delta _1$ is mapped into a disc in $S^2$, of radius $\leq \delta _2$. Then,
by Lemma \ref{lem7}, there is a positive number $\delta $ such that any
solid cylinder of radius $\leq \delta $ and based at a point in $K$, is
supported on two arcs of length $\leq \delta _1$.\hfill 
\endproof%

Recall that, if $X$ is a subset of $\mathcal{H}^2$, we denote by $G(X)$ the
set of geodesics in $\mathcal{H}^2$ which intersect $X$. To simplify
notation, we will write $G(x)$ for the set of geodesics through the point $%
x\in \mathcal{H}^2$, and $G(x,y)$ for the set of geodesics intersecting the
open geodesic segment $(x,y)$.

If $\Gamma $ is a group of isometries of $\mathcal{H}^2$, we denote by $%
G_{\Gamma}^{\prime }$ the set of geodesics in $\mathcal{H}^2$ which do not
intersect any of their translates by $\Gamma $: 
\[
G_{\Gamma }^{\prime }=\{\gamma \in G(\mathcal{H}^2):\forall g\in \Gamma
,g(\gamma )\cap \gamma =\emptyset \text{ or }g(\gamma )=\gamma \}. 
\]

In the following Lemma we consider the angle between unoriented geodesics to
lie in the interval $[0,\frac \pi 2]$.

\begin{lemma}
\label{lem9}Let $\ell $ and $\theta $ be positive numbers. Then there is a
positive number $\zeta $ with the following property. Let $x,y\in \mathcal{H}%
^2$, $\gamma $ the geodesic carrying the segment $[x,y]$, $g\in PSL(2,\mathbb{R)%
}$ and $\gamma ^{\prime }\in G_{\langle g\rangle }^{\prime }$, such that:%

\begin{itemize}
\item[\rm i)] The hyperbolic distance $d(x,y)\leq \ell $.
\item[\rm ii)] The geodesic segments $[x,y]$ and $[g(x),g(y)]$ intersect, and the angle
between $\gamma $ and $g(\gamma )$ is $\alpha \geq \theta $.
\item[\rm iii)] $\gamma ^{\prime }$ intersects the segment $[x,y]$ and the angle
between $\gamma $ and $\gamma ^{\prime }$ is $\beta $.
\end{itemize}
Then $\beta \geq \zeta $.
\end{lemma}

\proof
Without loss of generality, we may asume that $x=i\in \mathcal{H}^2$ and $%
y=ti$. The angle of intersection between the geodesics $\delta $ and $%
g(\delta )$ is a continuous function of $\delta $. Hence there is a
neighbourhood $U$ of $\gamma \in G(\mathcal{H}^2)$ disjoint from $G_{\langle
g\rangle }^{\prime }$, that is consisting of geodesics $\delta $ such that $%
g(\delta )$ intersects $\delta $.

There is a positive number $r$ such that the (two dimensional) solid
cylinder $C(\gamma ,i\sqrt{t},r)$ has the property: if $\delta \subset
C(\gamma ,i\sqrt{t},r)$ then $\delta \in U$. Then it is easy to show, using
hyperbolic trigonometry, that there is a positive number $\zeta $ such that
any geodesic $\delta $ intersecting $[x,y]$ at an angle $\leq \zeta $ is
contained in $C(\gamma ,i\sqrt{t},r)$, and hence $\delta \notin G_{\langle
g\rangle }^{\prime }$.\hfill 
\endproof%

\begin{corollary}
\label{cor10} If $g$ is a hyperbolic isometry of $\mathcal{H}^2$ and $x\in 
\mathcal{H}^2$ does not lie on the axis of $g$, then there is a positive
number $\zeta $ with the following property. If $\mu $ is any geodesic
lamination invariant by $g$, then no leaf of the lamination intersects the
geodesic segment $[x,g(x)]$ at an angle smaller than $\zeta $.\hfill 
\endproof%
\end{corollary}

\begin{lemma}
\label{lem11} Let $\ell ,\theta $ and $\varepsilon $ be positive numbers.
Then there is a positive number $r$ with the following property. Let $x,y\in 
\mathcal{H}^2$ with $d(x,y)\leq \ell \,$, and let $\gamma $ be the geodesic
carrying the segment $[x,y]$. Let $g\in PSL(2,\mathbb{R)}$ be such that $[x,y]$
intersects $[g(x),g(y)]$ at the point $x_0$, and at an angle $\alpha \geq
\theta $. If $\delta \in G_{\langle g\rangle }^{\prime }\cap G(D(x_0,r))$,
then $\delta $ intersects both $\gamma $ and $g(\gamma )$, and the points of
intersection lie in $D(x_0,\varepsilon )$.
\end{lemma}

\proof%
Since $g^{-1}(x_0)\in [x,y]$, we have $d(g^{-1}(x_0),x_0)\leq \ell $. We
consider the geodesic segment $[x^{\prime },y^{\prime }]$ of length $3\ell $
on the geodesic $\gamma $, centred at $x_0$.

Let $U$ be a neighbourhood of $\gamma \in G(\mathcal{H}^2)$ disjoint from $%
G_{\langle g\rangle }^{\prime }$. There is $r_1$ such that any geodesic
which intersects $D(x_0,r_1)$ and does not intersect $[x^{\prime },y^{\prime
}]$, lies in $U$, and hence it is not in $G_{\langle g\rangle }^{\prime }$.
So, if $\delta \in G_{\langle g\rangle }^{\prime }\cap G(D(x_0,r_1))$, $%
\delta $ intersects the segment $[x^{\prime },y^{\prime }]$. Similarly,
there is $r_2$ such that if $\delta \in G_{\langle g\rangle }^{\prime }\cap
G(D(x_0,r_2))$, $\delta $ intersects the segment $[g(x^{\prime
}),g(y^{\prime })]$.

By Lemma \ref{lem9}, the angle at the points of intersection is greater than
a constant $\zeta $. If $r$ satisfies $0<r<\min (r_1,r_2)$ and $\sinh r<\sin
\zeta \sinh \varepsilon $, then it has the required property.\hfill 
\endproof%

The following Lemma shows that, under certain conditions, taking integrals
along geodesic segments describes weak convergence of measures.

\begin{lemma}
\label{lem12}Let $\{\mu _n\}$ be a sequence of measured geodesic laminations
on $\mathcal{H}^2$, invariant by $g\in PSL(2,\mathbb{R)}$, and assume that $\mu
_n$ converge weakly to a measured lamination $\mu $. Let $\gamma $ be a
geodesic in $\mathcal{H}^2$, such that $\gamma $ and $g(\gamma )$ intersect
at one point. Then, for every geodesic segment $[u,v]$ on $\gamma $ and for
every continuous function $f\co [u,v]\rightarrow [0,1]$, with $f(u)=f(v)=0$,
the sequence $\int_{[u,v]}f\mu _n$ converges to $\int_{[u,v]}f\mu $.
\end{lemma}

\proof%
Since $\gamma $ intersects $g(\gamma )$ at one point, there is a
neighbourhood $U$ of $\gamma $ in $G(\mathcal{H}^2)$ which is disjoint from $%
G_{\langle g\rangle }^{\prime }$. We define a continuous function $\tilde{f}%
\co G(\mathcal{H}^2)\rightarrow [0,1]$ by letting $\tilde{f}(\delta )=f(y)$ if $%
y\in [u,v]$ and $\delta \in G(y)-U$, and extending continuously to the rest
of $G(\mathcal{H}^2)$. Then, for any measured geodesic lamination $\nu $
invariant by $g$, 
\[
\tilde{f}\nu (G(u,v))=\int_{[u,v]}f\nu.
\]
\vglue -.8cm
\endproof

\section{The theorems\label{sec3}}

We fix a reference point $[\rho _0]\in T(S)$, and we consider a point $[\rho
]\in Q(S)$. Let $g_1,\ldots ,g_k\in PSL(2,\mathbb{R)}$ be a set of generators
for $\Gamma _0=\rho _0(\pi _1(S))$. Let $x\in \mathcal{H}^2$ be a point
which does not lie on the axis of any conjugate of the generators $g_j$.

Let $\theta $ be the minimum of the angles between the geodesics carrying
the segments $[g_j^{-1}(x),x]$ and $[x,g_j(x)]$, for $j=1,\ldots ,k$. Let $d$
and $d^{\prime }$ be the maximum and the minimum, respectively, of the
distances between $x$ and $g_j(x)$, for $j=1,\ldots ,k$.

Let $K$ be a compact disc in $\mathcal{H}^2$ containing in its interior the
points $x$, $g_j(x)$, $g_j^{-1}(x)$, for $j=1,\ldots ,k$, and projecting
onto $S_0=\mathcal{H}^2/\Gamma _0$. Let $L$ be a compact disc in $\mathcal{H}%
^3$ such that $\phi _{*}(G(K))\subset G(L)$.

We consider a positive integer $m$, and a positive number $r(m)$ such that $%
d/m$ is less than the number $\delta (K,L,r(m))$ given by Lemma \ref{lem8}.

Let $\mu $ be a complex measured geodesic lamination on $\mathcal{H}^2$,
invariant by the group $\Gamma _0$, with $||\mu ||<M_0$. We consider one of
the generators $g_j$, $j=1,\ldots ,k$, and to simplify notation we drop the
suffix $j$ for the time being. Let $\gamma $ denote the geodesic carrying
the segment $[x,g(x)]$. We divide the segment $[x,g(x)]$ into $m$ equal
subsegments, by the points 
\[
x=x_0,x_1,\ldots ,x_{m-1},x_m=g(x). 
\]

If $[x,y]$ is a geodesic segment in $\mathcal{H}^2$ and $\nu $ is a measure
on a set of geodesics in $\mathcal{H}^2$, we introduce the notation 
\[
\int_{\lbrack x,y]}^{\prime }\nu =\frac 12\nu \left( G\left( x\right)
\right) +\nu \left( G\left( x,y\right) \right) +\frac 12\nu \left( G\left(
y\right) \right) 
\]

We define two new measures on the set $G(\mathcal{H}^2)$ of geodesics in $%
\mathcal{H}^2$ in the following way. For every $i=1,\ldots ,m$, let $\tilde{%
\gamma}_i$ be a geodesic in $\limfunc{supp}\mu $, intersecting $\gamma $ in $%
[x_{i-1},x_i]$. We define, for $i=1,\ldots ,m$, 
\[
\tilde{\mu }(\tilde{\gamma}_i)=\int_{[x_{i-1},x_i]}^{\prime }\mu. 
\]
For every $i=1,\ldots ,m-1$, let $\gamma _i^{\prime }$ be the geodesic in $%
\limfunc{supp}\mu $ intersecting the open segment $(x_{i-1},x_{i+1})$ as
near as possible to $x_i$. Let $\lambda _i\co [x_0,x_m]\rightarrow [0,1]$, $%
i=1,\ldots ,m-1$, be continuous functions satisfying

\begin{enumerate}
\item  $\limfunc{supp}(\lambda _i)\subset [x_{i-1},x_{i+1}]$ and

\item  $\sum_{i=1}^{m-1}\lambda _i\left( x\right) =1$ for all $x\in [x_0,x_m]
$.
\end{enumerate}

Then, in particular, $[x_0,x_1]\subset \lambda _i^{-1}(1)$ and $%
[x_{m-1},x_m]\subset \lambda _{m-1}^{-1}(1)$. We define, for $i=1,\ldots
,m-1 $, 
\[
\mu ^{\prime }(\gamma _i^{\prime })=\int_{[x_{i-1},x_{i+1}]}\lambda _i\mu 
\]

Now we define 
\[
C_i=A(\phi _{*}(\tilde{\gamma}_i),\tilde{\mu}(\tilde{\gamma}_i))\hspace{1cm}%
\text{for }i=1,\ldots ,m 
\]
and 
\[
D_i=A(\phi _{*}(\gamma _i^{\prime }),\mu ^{\prime }(\gamma _i^{\prime }))%
\text{\hspace{1cm}for }i=1,\ldots ,m-1. 
\]
We want to bound the norm $||C_1C_2\cdots C_m-D_1D_2\cdots D_{m-1}||$.

We put $a_i=\int_{[x_{i-1},x_i]}^{\prime }\lambda _i\mu $ and $%
b_i=\int_{[x_i,x_{i+1}]}^{\prime }\lambda _i\mu $. Then $\mu ^{\prime
}(\gamma _i^{\prime })=a_i+b_i$, for $i=1,\ldots ,m-1$, and $\tilde{\mu }(%
\tilde{\gamma}_1)=a_1$, $\tilde{\mu }(\tilde{\gamma}_m)=b_{m-1}$, and for $%
i=2,\ldots ,m-1$, $\tilde{\mu }(\tilde{\gamma}_i)=b_{i-1}+a_i$.

We put $D_i^l=A(\phi _{*}(\gamma _i^{\prime }),a_i)$ and $D_i^r=A(\phi
_{*}(\gamma _i^{\prime }),b_i)$. With this notation we have 
\begin{eqnarray*}
\lefteqn{\left\| C_1\cdots C_m-D_1\cdots D_{m-1}\right\| \leq }\hspace{1cm}
\\
&&\left\| C_1\cdots C_{m-1}\right\| \left\| C_m-D_{m-1}^r\right\| \\
&&\mbox{}+\left\| C_1\cdots C_{m-2}\right\| \left\|
C_{m-1}-D_{m-2}^rD_{m-1}^l\right\| \left\| D_{m-1}^r\right\| \\
&&\mbox{}+\cdots +\left\| C_1\cdots C_{s-1}\right\| \left\|
C_s-D_{s-1}^rD_s^l\right\| \left\| D_s^rD_{s+1}\cdots D_{m-1}\right\| \\
&&\mbox{}+\cdots +\left\| C_1-D_1^l\right\| \left\| D_1^rD_2\cdots
D_{m-1}\right\|.
\end{eqnarray*}
Then, by Lemma \ref{lem1.1}, there is a positive number $M_1$, depending on $%
L$ and $M_0$, which is an upper bound for the norm of the factors of the
form $C_1\cdots C_s$, $D_s^rD_{s+1}\cdots D_{m-1}$. By Lemma \ref{lem4.1},
there is a positive number $M_2$, depending on $L$ and $M_0$, such that each
factor of the form $C_s-D_{s-1}^rD_s^l$ has norm bounded by $M_2r(m)\tilde{%
\mu}(\tilde{\gamma}_s)$. Then 
\begin{eqnarray}
\left\| C_1\cdots C_m-D_1\cdots D_{m-1}\right\| & \leq & M_0M_1^2M_2r(m) 
.  \label{eq3}
\end{eqnarray}

In the following we want to examine the behaviour of $D_1\cdots D_{m-1}$ as $%
m\rightarrow \infty $ and as the lamination $\mu $ changes. For this we must
consider more carefully the leaves of the lamination near $x$.

By Lemma \ref{lem11}, there is an open set $U\subset G(K)$, depending on $%
d,\theta $ and $d^{\prime }/m$ such that, if $\delta $ is any geodesic in $%
U\cap \limfunc{supp}\mu $, then $\delta $ intersects the geodesics $\gamma $
and $g(\gamma )$ at a distance less than $d^{\prime }/m$ from $x$. Let $\chi
\co G(K)\rightarrow [0,1]$ be a continuous function, with $\limfunc{supp}\chi
\subset U$ and $\chi |_{G(x)}=1$. We introduce the notation 
\[
a^{\prime }=\int_{[x_0,x_1]}\chi \mu \hspace{2cm}\bigskip a^{\prime \prime
}=\int_{[x_0,x_1]}^{\prime }(1-\chi )\mu 
\]
\[
b^{\prime }=\int_{[x_{m-1},x_m]}(\chi \circ g^{-1})\mu \hspace{2cm}\bigskip
b^{\prime \prime }=\int_{[x_{m-1},x_m]}^{\prime }(1-\chi \circ g^{-1})\mu 
\]
\[
P=A(\phi _{*}(\gamma _1^{\prime }),a^{\prime })\hspace{2cm}\bigskip Q=A(\phi
_{*}(\gamma _1^{\prime }),a^{\prime \prime }) 
\]
\[
R=A(\phi _{*}(\gamma _{m-1}^{\prime }),b^{\prime \prime })\hspace{2cm}%
\bigskip S=A(\phi _{*}(\gamma _{m-1}^{\prime }),b^{\prime }), 
\]
and we have 
\[
D_1=PQD_1^r\hspace{2cm}\bigskip D_{m-1}=D_{m-1}^lRS. 
\]

Let $\{\mu _n\}$ be a sequence of complex measured geodesic laminations on
the surface $S_0$, converging weakly in $\mathcal{M}\left( G\left( K\right)
\right) $ to a measured lamination $\mu _0$. Then, by the Uniform
Boundedness Principle, there is a positive number $M_0$ such that $||\mu
_n||\leq M_0$ for all $n\geq 0$.

For each positive integer $m$, for each $i=1,\ldots ,m-1$, for each $%
j=1,\ldots ,k$ and for each measured lamination $\mu _n$, $n\geq 0$, we
define as above the points $x_{j,m,i}$, the geodesics $\gamma
_{n,j,m,i}^{\prime }$, the functions $\lambda _{j,m,i}$, the quantities $%
a_{n,,j,m,i}$, $b_{n,j,m,i}$, $a_{n,j,m}^{\prime }$, $b_{n,j,m}^{\prime }$
and the isometries $D_{n,j,m,i}$, $P_{n,j,m}$, $Q_{n,j,m}$, $R_{n,j,m}$, $%
S_{n,j,m}$.

Let $B_{n,j,m}=D_{n,j,m,1}\cdots D_{n,j,m,m-1}$. We want to find a bound for
the norm of the difference between $B_{0,j,m}g_j$ and some conjugate of $%
B_{n,j,m}g_j$.

\begin{lemma}
\label{lemmain}With the above notation, there exist positive numbers $N_1,N_2
$ and functions $r\co \mathbb{N}\rightarrow \mathbb{R}$, $\varepsilon \co \mathbb{N\times N}%
\rightarrow \mathbb{R}$ such that 
\[
\lim_{m\rightarrow \infty }r(m)=0,\qquad \lim_{n\rightarrow \infty
}\varepsilon (m,n)=0\text{\quad for each }m\in \mathbb{N}
\]
and 
\[
\left\|
P_{0,1,m}P_{n,1,m}^{-1}B_{n,j,m}g_jP_{n,1,m}P_{0,1,m}^{-1}-B_{0,j,m}g_j%
\right\| \leq N_1r(m)+N_2\varepsilon (m,n).
\]
\end{lemma}

\proof%
To simplify notation, we drop the index $m$ for the time being, and write,
for example, $D_{n,j;i}$ for $D_{n,j,m,i}$. We have 
\begin{eqnarray}
\lefteqn{\left\| P_{0,1}P_{n,1}^{-1}B_{n,j}g_jP_{n,1}P_{0,1}^{-1}-B_{0,j}
g_j\right\| \leq } \hspace{1cm}  \nonumber \\
&&\left\|P_{0,1}P_{n,1}^{-1}B_{n,j}g_jP_{n,1}P_{0,1}^{-1}-P_{0,j}
P_{n,j}^{-1}B_{n,j}g_jP_{n,j}P_{0,j}^{-1}\right\|  \label{eq1} \\
&&\mbox{}+\left\|
P_{0,j}P_{n,j}^{-1}B_{n,j}g_jP_{n,j}P_{0,j}^{-1}g_j^{-1}-P_{0,j}
P_{n,j}^{-1}B_{n,j}S_{n,j}^{-1}S_{0,j}\right\| \left\| g_j\right\|  \nonumber
\\
&&\mbox{}+\left\|
P_{0,j}P_{n,j}^{-1}B_{n,j}S_{n,j}^{-1}S_{0,j}-B_{0,j}\right\| \left\|
g_j\right\| .  \nonumber
\end{eqnarray}
We will find upper bounds for the three terms of the right hand side of the
above inequality.

The first term of (\ref{eq1}) is bounded above by 
\begin{eqnarray}\lefteqn{
\left\| P_{0,1}P_{n,1}^{-1}-P_{0,j}P_{n,j}^{-1}\right\| \left\|
B_{n,j}g_jP_{n,1}P_{0,1}^{-1}\right\|}\hspace{1cm}  \nonumber \\
&&
+\left\|
P_{0,j}P_{n,j}^{-1}B_{n,j}g_j\right\| \left\|
P_{n,j}P_{0,j}^{-1}-P_{n,j}P_{0,j}^{-1}\right\|. 
\nonumber
\end{eqnarray}
By Lemma \ref{lem1.1}, the factors containing $g_j$ are bounded above by $%
M_1 $. We consider the other factor in each term. Recall that $%
P_{n,j}=A(\phi _{*}(\gamma _{n,j;1}^{\prime }),a_{n,j}^{\prime })$. We have 
\begin{eqnarray}
\lefteqn{\left\| P_{0,j}P_{n,j}^{-1}-P_{0,1}P_{n,1}^{-1}\right\| \leq} 
\hspace{1cm}  \nonumber \\
&&\left\| P_{0,j}\right\| \left\| P_{n,j}^{-1}-A(\phi _{*}(\gamma
_{0,j;1}^{\prime }),-a_{n,j}^{\prime })\right\|  \label{eq2} \\
&&\mbox{}+\left\| A(\phi _{*}(\gamma _{0,j;1}^{\prime }),a_{0,j}^{\prime
}-a_{n,j}^{\prime })-A(\phi _{*}(\gamma _{0,1;1}^{\prime }),a_{0,1}^{\prime
}-a_{n,1}^{\prime })\right\|  \nonumber \\
&&\mbox{}+\left\| P_{0,1}\right\| \left\| A(\phi _{*}(\gamma
_{0,1;1}^{\prime }),-a_{n,1}^{\prime })-P_{n,1}^{-1}\right\| . 
\nonumber
\end{eqnarray}
By Lemma \ref{lem4}, there is a positive constant $M^{\prime }$ such that
the first and the third term of the right hand side of (\ref{eq2}) are
bounded by $M_0M_1M^{\prime }r(m)$. To find a bound for the second term we
consider two cases.

\begin{enumerate}
\item  The segment $[x_0,x_{j;1}]$ intersects the same geodesics in $%
\limfunc{supp}(\chi \mu _n)$ as does the segment $[x_0,x_{1;1}]$.

\item  The two segments intersect different sets of geod\-esics in $\limfunc{%
supp}(\chi \mu _n)$.
\end{enumerate}

Let $z_{n,i}=\int_{[x_0,x_{i;1}]}\chi (\mu _0-\mu _n)=a_{0,i}^{\prime
}-a_{n,i}^{\prime }$.

In case (1), $z_{n,j}=z_{n,1}$, and the geodesics $\gamma _{0,j;1}^{\prime
},\gamma _{0,1;1}^{\prime }$ lie in a (2--dimensional) solid cylinder of
radius $d/m$ based at $x_0$. The segments $[x_0,x_{j;1}]$ and $[x_0,x_{1;1}]$
induce concurrent orientations on the geodesics $\gamma _{0,j;1}^{\prime }$
and $\gamma _{0,1;1}^{\prime }$ respectively. So, by Lemma \ref{lem4}, 
\[
\left\| A(\phi _{*}(\gamma _{0,j;1}^{\prime }),z_{n,j})-A(\phi _{*}(\gamma
_{0,1;1}^{\prime }),z_{n,1}\right\| \leq M_0M^{\prime }r(m). 
\]

Note that if $\mu _n$ satisfies the conditions of case (1) for large enough $%
n$, then $\mu _0$ also satisfies these conditions.

In case (2), the orientations induced by the segments $[x_0,x_{j;1}]$ and $%
[x_0,x_{1,1}]$ on the geodesics $\gamma _{0,j;1}^{\prime }$ and $\gamma
_{0,1;1}^{\prime }$ respectively, are not concurrent. Hence, by Lemma \ref
{lem4}, 
\[
\left\| A(\phi _{*}(\gamma _{0,j;1}^{\prime }),z_{n,j})-A(\phi _{*}(\gamma
_{0,1;1}^{\prime }),z_{n,1)}\right\| \leq M_0M^{\prime }r(m)+M^{\prime
\prime }|z_{n,j}+z_{n,1}|. 
\]
Note that, in this case, 
\[
a_{0,j}^{\prime }+a_{0,1}^{\prime }=\int_{[x_0,x_{j;1}]}\chi \mu
_0+\int_{[x_0,x_{1;1}]}\chi \mu _0=\chi \mu _0(G) 
\]
and similarly for $\mu _n$. Hence $z_{n,j}+z_{n,1}=\chi \mu _0(G)-\chi \mu
_n(G)$. Let 
\[
\varepsilon _0(m,n)=\sup_{s\geq n}|\chi _m\mu _0(G)-\chi _m\mu _s(G)|. 
\]

Now we turn our attention to the second term of equation (\ref{eq1}). This term
involves only the generator $g_j$, so we drop the subscript $j$ from the
notation. We have 
\begin{eqnarray*}
\lefteqn{\left\|
P_0P_n^{-1}B_ngP_nP_0^{-1}g^{-1}-P_0P_n^{-1}B_nS_n^{-1}S_0\right\| \leq } \\
&&\left\| P_0P_n^{-1}B_n\right\| \left\| S_n^{-1}\right\| \left\|
S_ngP_n^{-1}g^{-1}-S_0gP_0g^{-1}\right\| \left\| gP_0^{-1}g^{-1}\right\| 
.
\end{eqnarray*}
We consider the term $S_ngP_n^{-1}g^{-1}$, which is equal to 
\[
A\left( \phi _{*}(\gamma _{n;m-1}^{\prime }),\int_{[x_{;m-1},x_{;m}]}(\chi
\circ g^{-1})\mu _n\right) A\left( \phi _{*}(g(\gamma _{n;1}^{\prime
}),\int_{[x_0,x_{;1}]}\chi \mu _n\right). 
\]
Since $\mu _n$ is invariant by $g$, and $x_{;m}=g(x_0)$, we have 
\[
\int_{\lbrack x_{;m},g(x_{;1})]}(\chi \circ g^{-1})\mu
_n=\int_{[x_0,x_{;1}]}\chi \mu _n. 
\]

We have to consider two cases:

\begin{enumerate}
\item  The segments $[x_{;m-1},x_{;m}]$ and $[x_{;m},g(x_{;1})]$ intersect
the same geodesics in $\limfunc{supp}((\chi \circ g^{-1})\mu _n)$.

\item  The segments $[x_{;m-1},x_{;m}]$ and $[x_{;m},g(x_{;1})]$ intersect
different sets of geod\-esics in $\limfunc{supp}((\chi \circ g^{-1})\mu _n)$.
\end{enumerate}

In case (1), we let $z_n=\int_{[x_{;m-1},x_{;m}]}(\chi \circ g^{-1})\mu
_n=\int_{[x_{;m},g(x_{;1})]}(\chi \circ g^{-1})\mu _n$. The geodesics $%
\gamma _{n;m-1}^{\prime }$and $g(\gamma _{n;1}^{\prime })$lie in a solid
cylinder of radius $d/m$, based at $x_{;m}$, and the orientations induced by
the segments $[x_{;m-1},x_{;m}]$ and $[x_{;m},g(x_{;1})]$ are not
concurrent. Hence, by Lemma \ref{lem4.1}, $\left\| S_ngP_ng^{-1}-I\right\|
\leq M_0M_2r(m)$. As before, if $\mu _n$ satisfies the conditions of case
(1) for large enough $n$, then $\mu _0$ also satisfies these conditions.
Hence 
\[
\left\| S_ngP_ng^{-1}-S_0gP_0g^{-1}\right\| \leq 2M_0M_2r(m). 
\]

In case (2), since $\mu _n$ is invariant by $g$, and $x_{;m}=g(x_0)$, we
have 
\[
\int_{\lbrack x_{;m},g(x_{;1})]}(\chi \circ g^{-1})\mu
_n+\int_{[x_{;m-1},x_{;m}]}(\chi \circ g^{-1})\mu _n=\chi \mu _n(G) 
\]
and if $n$ is large enough, the same is true of $\mu _0$. Then 
\begin{eqnarray*}
\lefteqn{\left\| S_ngP_ng^{-1}-S_0gP_0g^{-1}\right\| \leq } \hspace{1cm} \\
&&\left\| S_ngP_ng^{-1}-A(\phi _{*}(\gamma _{n;m-1}^{\prime }),\chi \mu
_n(G))\right\| \\
&&\mbox{}+\left\| A(\phi _{*}(\gamma _{n;m-1}^{\prime }),\chi \mu
_n(G))-A(\phi _{*}(\gamma _{0;m-1}^{\prime }),\chi \mu _0(G))\right\| \\
&&\mbox{}+\left\| A(\phi _{*}(\gamma _{0;m-1}^{\prime }),\chi \mu
_0(G))-S_0gP_0g^{-1}\right\|.
\end{eqnarray*}
By Lemma \ref{lem4} and Lemma \ref{lem4.1}, this is bounded above by $%
M^{\prime }r(m)+M^{\prime \prime }\varepsilon (m,n)$.

The third term of equation (\ref{eq1}) is bounded by 
\[
\left\| P_0\right\| \left\| P_n^{-1}B_nS_n^{-1}-P_0^{-1}B_0S_0^{-1}\right\|
\left\| S_0\right\| \left\| g\right\|.
\]
But 
\begin{eqnarray*}
\lefteqn{\left\| P_n^{-1}B_nS_n^{-1}-P_0^{-1}B_0S_0^{-1}\right\| =} \\
&&\left\| Q_nD_{n;1}^rD_{n;2}\cdots
D_{n;m-2}D_{n;m-1}^lR_n-Q_0D_{0;1}^rD_{0;2}\cdots
D_{0;m-2}D_{0;m-1}^lR_0\right\|
\end{eqnarray*}
and by Lemma \ref{lem1.1}, this is bounded by 
\begin{eqnarray}
&&M_1^2\mbox{\huge (} \left\| D_{n;m-1}^lR_n-D_{0;m-1}^lR_0\right\|
+\sum_{i=2}^{m-2}\left\| D_{n,i}-D_{0,i}\right\| +  \nonumber \\
&&\hspace{2cm} \mbox{}+\left\| Q_nD_{n;1}^r-Q_0D_{0;1}^r\right\| 
\mbox{\huge
)}.  \label{eq4}
\end{eqnarray}
Note that $Q_nD_{n;1}^r=A\left( \phi _{*}(\gamma _{n;1}^{\prime
}),\int_{[x_0,x_{;1}]}\lambda _{;1}(1-\chi )\mu _n\right) $ and hence 
\[
\left\| Q_nD_{n;1}^r-Q_0D_{0;1}^r\right\| \leq M^{\prime }r(m)+M^{\prime
\prime }\varepsilon _1(m,n) 
\]
where $\varepsilon _1(m,n)=\sup_{s\geq n}\left| \int_{[x_0,x_{;1}]}\lambda
_{;1}(1-\chi _m)(\mu _s-\mu _0)\right| $, and similarly for the other terms
of (\ref{eq4}), for suitable $\varepsilon _i$, $i=2,\ldots ,m-1$.

To complete the proof of Lemma \ref{lemmain} we must show that $r(m)$ and $%
\varepsilon (m,n)=\sum_{i=0}^{m-1}\varepsilon _i(m,n)$ have the required
properties. It is clear that we can choose a sequence $r(m)$, with $%
\lim_{m\rightarrow \infty }r(m)=0$, such that the pair $r=r(m)$, $\delta
=d/m $ satisfy the conditions of Lemma \ref{lem8}. Lemma \ref{lem12} implies
that, for each $m$, $\lim_{n\rightarrow \infty }\varepsilon (m,n)=0$. \hfill%
\endproof%

We let $E_{n,j,m}=C_{n,j,m,1}\cdots C_{n,j,m,m}$ and $%
H_{n,m}=P_{0,1,m}P_{n,1,m}^{-1}$. Then, combining the above result with (\ref
{eq3}), we have 
\begin{equation}
\left\| H_{n,m}E_{n,j,m}g_jH_{n,m}^{-1}-E_{0,j,m}g_j\right\| \leq
M(r(m)+\varepsilon (m,n).  \label{eq5}
\end{equation}
If $g_1,\ldots ,g_k$ is a set of generators for $\Gamma _0$, the space $R$
of homomorphisms $\rho \co \Gamma _0\rightarrow PSL(2,\mathbb{C)}$ with
quasi-Fuchsian image is a subspace of $PSL(2,\mathbb{C)}^k$, and $Q(S)$ is a
subspace of the quotient by the adjoint action on the left, $\left. PSL(2,%
\mathbb{C)}^k\right/ PSL(2,\mathbb{C)}$. Let 
\[
\rho _{n,m}=\left( H_{n,m}E_{n,j,m}g_jH_{n,m}^{-1},\quad j=1,\ldots
,k\right) 
\]
\[
\rho _{n,m}=\left( E_{0,j,m}g_j,\quad j=1,\ldots ,k\right) 
\]
and let $[\rho _{n,m}]$ denote the equivalence class of $\rho _{n,m}$ in $%
\left. PSL(2,\mathbb{C)}^k\right/ PSL(2,\mathbb{C)}$.

Let $n(m)$ be a sequence such that $n(m)\geq m$ and $\varepsilon
(n(m),m)\leq 1/m$. Then $\lim_{m\rightarrow \infty }\rho _{n(m),m}=\rho
_{\mu _0}$. As $m\rightarrow \infty $, $[\rho _{n,m}]$ converge, uniformly
in $n$, to the bending deformation $[\rho _{\mu _n}]$, \cite{K2}. Hence, $%
\lim_{m\rightarrow \infty }[\rho _{n(m),m}]=\lim_{m\rightarrow \infty }[\rho
_{\mu _{n(m)}}]=\lim_{n\rightarrow \infty }[\rho _{\mu _n}]$, and we have 
\begin{equation}
\lim_{n\rightarrow \infty }[\rho _{\mu _n}]=[\rho _{\mu _0}].
\label{eq6}
\end{equation}

To complete the proof of Theorem \ref{thm1}, it remains to show that the
convergence is uniform in compact subsets of $\mathcal{D}$. If $([\rho
],t)\in \mathcal{D}$, each bound used in the proof of (\ref{eq6}) depends at
most linearly on $t$, while it depends on $\rho $ only in terms of the
endpoints of a finite number of geodesics $\phi _{*}(\gamma )$. The
endpoints of the geodesic $\phi _{*}(\gamma )$ are, for each $\gamma $,
holomorphic functions of $[\rho ]$. Hence each bound can be chosen uniformly
on each compact subset of $\mathcal{D}$.

Note that $\mathcal{D}$ contains in its interior the set $Q(S)\times \{0\}$.
If the laminations $\mu _n$ are real for all but a finite number of $n$,
then $\mathcal{D}$ also contains the set $Q(S)\times \mathbf{R}$, but this
is not true in the general case.

To prove Theorem \ref{thm2} we recall that the bending vector field $T_{\mu }$
is defined by 
\[
T_{\mu }([\rho ])=\frac \partial {\partial t}B_{\mu }([\rho ],t).
\]
The vector fields $T_{\mu _n}$ are holomorphic, and $B_{\mu _n}([\rho ],t)$
converge to $B_{\mu _0}([\rho ],t)$ for $([\rho ],t)\in \mathcal{D}$. It
follows that $T_{\mu _n}$ converge to $T_{\mu _0}$, uniformly on compact
subsets of $Q(S)$.

We conclude with the proof of Theorem \ref{thm3}. We consider the subset of $%
\mathcal{ML}(S)$ consisting of measured laminations with non negative real
and imaginary parts, and we denote it by $\mathcal{ML}^{++}(S)$. We identify 
$\mathcal{ML}^{++}(S)$ with a subset of the set of pairs of positive
measured laminations $\mathcal{ML}_{\mathbb{R}}^{+}(S)\times \mathcal{ML}_{\mathbb{%
R}}^{+}(S)$. If $\nu \in \mathcal{ML}^{++}(S)$, then $\func{Re}\nu $ and $%
\func{Im}\nu $ are in $\mathcal{ML}_{\mathbb{R}}^{+}(S)$ and they satisfy the
condition 
\begin{equation}
\limfunc{supp}(\func{Re}\nu )\cup \limfunc{supp}(\func{Im}\nu )\text{ is a
geodesic lamination.}  \label{eq7}
\end{equation}
Conversely, any pair $\nu _1,\nu _2$ of positive measured laminations
satisfying (\ref{eq7}) define a measure $\nu =\nu _1+i\nu _2\in \mathcal{ML}%
^{++}(S)$. The mapping is a homeomorphism of $\mathcal{ML}^{++}(S)$ onto a
subset of $\mathcal{ML}_{\mathbb{R}}^{+}(S)\times \mathcal{ML}_{\mathbb{R}}^{+}(S)$%
. But $\mathcal{ML}_{\mathbb{R}}^{+}(S)$ is homeomorphic to $\mathbb{R}^{6g-6}$, 
\cite{PH}. Thus $\mathcal{ML}^{++}(S)$ is first countable, and Theorem \ref
{thm2} implies that $\mu \mapsto T_{\mu }$ is continuous. Theorem \ref{thm3}
then follows by the continuity of the evaluation map.

\Addresses\recd
\end{document}